\documentclass[psamsfont,a4paper,12pt]{amsart}
\usepackage[english]{babel}   %loads english
\numberwithin{equation}{section}

\usepackage{setspace}     %controls line-spacing
\usepackage{enumerate} 	%allows lists
\usepackage{afterpage} 	%keeps tables and figures where they belong
\usepackage[usenames,dvipsnames]{color}
\usepackage{graphicx}	%allows inclusion of pictures

\usepackage{amsmath} 	%just apparently the most important thing ever
\usepackage{amssymb}	%symbols in math
\usepackage{mathrsfs} 	%really fancy script font
\usepackage{amsfonts} 	%standard fonts
\usepackage{latexsym} 	%more symbols
\usepackage[latin1]{inputenc} %sets up proper special characters

\usepackage{amsthm} 	%sets up theorem styling
\usepackage{stackrel}       %allows stacking in mathmode  
\usepackage{amscd} 	%commutative diagrams                                 

\usepackage{tabularx}	%automatically widening tables
\usepackage{longtable}	%for tables longer than a page

\usepackage{verbatim} %comment environment?
\usepackage{blindtext} %something about dummy text

\usepackage{amssymb,latexsym,amsmath}    % Standard packages
\usepackage{pstricks,multido,pst-plot}		 % Packages for drawing figures
\usepackage{multicol,fancyheadings}			   % Misc.
\usepackage{float}                       % Used to define Algorithm enviro
\usepackage[vcentermath,enableskew]{youngtab}

\usepackage{subfig}

\allowdisplaybreaks

\newcommand{\Z}{\mathbb{Z}} 
\newcommand{\Q}{\mathbb{Q}}

\newtheorem{theorem}{Theorem}[section]  %%with section numbering

\newtheorem*{conjecture}{Conjecture}

\newtheorem*{remark}{Remark}

\theoremstyle{definition}

\theoremstyle{remark}

\newtheorem{ind}[]{{\rm\it Indice}}

\usepackage{multicol}

\usepackage{tikz}
\usetikzlibrary{arrows}

\usepackage{paralist}
\usepackage{cite}

\title{Multiquadratic fields generated by characters of $A_n$}

\author[Dawsey]{Madeline Locus Dawsey*}
\address{Department of Mathematics,
Emory University, Atlanta, GA 30322}
\email{madeline.locus@emory.edu}

\author[Ono]{Ken Ono}
\address{Department of Mathematics, Emory University, Atlanta, GA 30322}
\email{ken.ono@emory.edu}

\author[Wagner]{Ian Wagner}
\address{Department of Mathematics, Emory University, Atlanta, GA 30332}
\email{iwagner@emory.edu}

\begin{document}

\thanks{*This author was previously known as Madeline Locus. The second author thanks the generous
support of the Asa Griggs Candler Fund, and the NSF (DMS-1601306).}

\begin{abstract} 
For a finite group $G$, let $K(G)$ denote the field generated over $\Q$ by its character values. For $n>24$, G. R. Robinson and J. G. Thompson \cite{RT} proved that $$K(A_n)=\Q\left (\{ \sqrt{p^*} \ : \ p\leq n \ {\text {\rm an odd prime with\ } p\neq n-2}\}\right),$$ where $p^*:=(-1)^{\frac{p-1}{2}}p$. Confirming a speculation of Thompson \cite{T}, we show that arbitrary suitable multiquadratic fields are  similarly generated by the values of $A_n$-characters restricted to elements whose orders  are only divisible by ramified primes. To be more precise, we say that a $\pi$-number is a positive integer whose prime factors belong to a set of odd primes $\pi:= \{p_1, p_2,\dots, p_t\}$. Let $K_{\pi}(A_n)$ be the field generated by the values of $A_n$-characters for even permutations whose orders are $\pi$-numbers. If $t\geq 2$, then we determine a constant $N_{\pi}$ with the property that for all $n> N_{\pi}$, we have $$K_{\pi}(A_n)=\Q\left(\sqrt{p_1^*}, \sqrt{p_2^*},\dots, \sqrt{p_t^*}\right).$$
\end{abstract}

\maketitle

\section{Introduction and Statement of Results}

If $G$ is a finite group, then let $K(G)$ denote the  field generated over $\Q$
by all of the $G$-character values. In stark contrast to the case of symmetric groups, where
$K(S_n)=\Q$, G. R. Robinson and J. G. Thompson \cite{RT} proved for alternating groups that
the $K(A_n)$ are generally large multiquadratic extensions.
In particular, for $n>24$, they proved that
$$
K(A_n)=\Q\left(\{ \sqrt{p^*} \ : \ p\leq n \ {\text {\rm an odd prime with\ } p\neq n-2}\}\right),
$$
where $m^*:=(-1)^{\frac{m-1}{2}}m$ for any odd integer $m$. 

In a letter \cite{T} to the second author in 1994, Thompson asked for a refinement of this result
that mirrors the  Kronecker-Weber Theorem and the theory of complex multiplication, where
abelian extensions are generated by the values of $e^{2\pi i x}$ and modular functions respectively, at arguments that determine the ramified primes.  Instead of employing special analytic functions at designated arguments, which is the gist of {\it Hilbert's 12th Problem}
\cite{H}, 
Thompson offered the characters of $A_n$ evaluated at elements whose orders are only divisible
by ramified primes.

To formulate this problem, we let $\pi:=\{p_1, p_2,\dots, p_t\}$ denote a set of $t\geq 2$  distinct odd primes\footnote{The phenomenon 
cannot hold
when $t=1$ for any $n\not\equiv 0, 1\pmod{p}$.} listed in increasing order. A
 $\pi$-{\it number} is a positive integer whose prime factors
belong to $\pi$. 
The speculation is that $K_{\pi}(A_n),$  the field generated by the values of $A_n$-characters restricted to
 elements $\sigma\in A_n$ with $\pi$-number order, generally generates $\Q(\sqrt{p_1^*}, \sqrt{p_2^*},\dots, \sqrt{p_t^*})$.
However, an inspection of the character tables for the first few $A_n$ casts doubt on this speculation. 
For example,  one easily finds that
$$
\Q=K_{\{3,5\}}(A_{22}) = K_{\{5,7\}}(A_{1738}) = K_{\{7,11\}}(A_{159557}) \ \ \ \ \ \ 
$$
by checking that  the $A_n$-character values for even permutations with $\pi$-number order are in $\Z$.
%%Note: There are no partitions of $n$ into distinct sums of these pi-numbers in these cases. Therefore Theorem 2.1
%%indicates that their character values at permutations with pi-number order are ordinary integers.

Despite these discouraging examples, Thompson's speculation is indeed true for sufficiently large $A_n$.
We let $\Omega_{\pi}:=p_t^{4\left(p_{t-1}\right)^2}$, and in turn we define
\begin{equation}\label{Npi1}
\mathcal{N}_{\pi}^{+}:= p_t + 2^{2^{2^{\Omega_{\pi}}}}.
\end{equation}
The formula for $\mathcal{N}_{\pi}^{+}$ provides a bound for this phenomenon.
It  can generally be improved when $t\geq 3$.
To make this precise,
we let $\mathbb{N}_{\pi}^{\square}=\{a_1, a_2,\dots\}$ be the subset of $\pi$-numbers that are perfect squares in increasing order, and we let
$S_{\pi}(n)$ denote the number of elements in $\mathbb{N}_{\pi}^{\square}$ not exceeding $n$.
  Choose a positive real number
$B_{\pi}\geq 2^{2^{40}}$ for which
$$
S_{\pi}(2n)-S_{\pi}(n) > 2^{20} \log_2 \log_2 n
$$
for all $n\geq B_{\pi}$. Then choose a positive real number $C_{\pi} \geq B_{\pi}$ for which
$$
\sum_{s=1}^{S_{\pi}(C_{\pi})} \| a_s \theta \|^2
 > 2 \log_2 B_{\pi} +50
$$
for all $\frac{1}{4B_{\pi}} \leq \theta \leq 1-\frac{1}{4B_{\pi}}$, where  $\| x \|$ is the distance between $x$ and the nearest
integer. In terms of $C_{\pi}$, we can choose
\begin{equation}\label{Npi2}
\mathcal{N}_{\pi}^{-}:=p_t+ \max \left \{C_{\pi}^{10}, 2^{2^{2^{400}}}\right \}.
\end{equation}
We obtain the following confirmation
of Thompson's speculation in terms of the bound $N_\pi:= \min\{ \mathcal{N}_{\pi}^{+}, \mathcal{N}_{\pi}^{-}\}$.

\begin{theorem}\label{Thm1}
Assuming the notation above, if $n> N_{\pi}$, then
$$
K_{\pi}(A_n)=\Q\left(\sqrt{p_1^*}, \sqrt{p_2^*},\dots, \sqrt{p_t^*}\right).
$$
\end{theorem}

\begin{remark}
The proof of Theorem~\ref{Thm1} uses the fact that 
 $n-p_i$, for each $1\leq i\leq t$, is the sum of distinct
squares of $\pi$-numbers provided $n > N_{\pi}$.  
These Diophantine conditions guarantee that  
$\Q(\sqrt{p_i^*}) \subset K_{\pi}(A_n)$. However, these conditions are not necessary for these inclusions.
\end{remark}

\begin{remark}
Although Theorem~\ref{Thm1} requires that $t\geq 2$, it is simple to generate $\Q(\sqrt{p^*})$ for an odd 
prime $p$ using permutations with cycle lengths that are $\pi$-numbers with $\pi=\{p,q\}$, where $q\neq p$ is an odd prime.
As mentioned above, if $n> N_{\pi}$, then $n-p$ is a sum of distinct squares of $\pi$-numbers. Then $K_{\pi_p}(A_n)$,
the field generated by the values of the $A_n$-characters restricted to permutations with a single cycle of length $p$ and other cycle lengths of such odd squares, satisfies $K_{\pi_p}(A_n)=\Q(\sqrt{p^*})$.
\end{remark}

To prove Theorem~\ref{Thm1}, we follow the straightforward approach of  Robinson and Thompson in \cite{RT}.
In Section 2.1, we recall standard facts about the characters of $A_n$. 
Theorem~\ref{Thm2Point1} is the key device
for relating cycle types to the surds $\sqrt{p^*}$. 
Then, in Section 2.2, we recall
a classical result of J. W. Cassels \cite{C} on partitions. We also recall
 recent work by J.-H. Fang and
Y.-G. Chen \cite{FC} on a problem of Erd\"os
 which implies that every 
large positive integer is the sum of distinct squares of $\pi$-numbers when $\pi=\{p_1, p_2\}$.
Theorem~\ref{Thm1} follows immediately from these results.

\section*{Acknowledgements} \noindent
The second author thanks John Thompson for sharing his speculation 25 years ago.
The authors thank Geoff Robinson for comments on the first draft of this paper.

\section{Nuts and Bolts and the proof of Theorem~\ref{Thm1}}

\subsection{Character theory for $A_n$}

It is well known that the representation theory of $S_{n}$ and $A_{n}$ can be completely described using the partitions of $n$.  In particular, a permutation $\sigma \in S_{n}$ has a cycle type that can be viewed as a partition of $n$, say $\lambda = (\lambda_{1}, \lambda_{2}, \dots, \lambda_{k})$.  The cycle type determines the conjugacy class of the permutation, and so the irreducible representations of $S_{n}$ can be indexed by (and constructed from) the partitions of $n$.  It is also well known that the only conjugacy classes which split in $A_{n}$ are those corresponding to partitions into distinct odd parts.  There is a bijection between the set of partitions $\lambda$ of $n$ into distinct odd parts and the set of self-conjugate partitions $\gamma$ of $n$ which is realized by identifying the parts of each $\lambda$ with the main hook lengths of some $\gamma$.  Theorem 2.5.13 of  \cite{JK} characterizes those $A_n$-character values that are not $\Z$-integral.

\begin{theorem}\label{Thm2Point1}
Let $\sigma \in A_n$ be a permutation with cycle type given by a partition $\lambda:=(\lambda_1,\lambda_2,\dots, \lambda_k)$ of $n$, and let $d_\lambda:=\prod_{i=1}^k\lambda_i$.
\begin{enumerate}
\item If $\lambda$ does not have distinct odd parts, then the $A_n$-character values of $\sigma$ are all in $\Z$.
\item If $\lambda$ has distinct odd parts, then let $\gamma$ be the self-conjugate partition of $n$ with main hook lengths $\lambda_1,\lambda_2,\dots,\lambda_k$, and let $\chi_\gamma$ be the $A_n$-character associated to $\gamma$.  We have that 
\begin{equation*}
\chi_{\gamma}(\sigma) = \frac{1}{2} \left( (-1)^{\frac{d_{\lambda}-1}{2}} \pm \sqrt{d_{\lambda}^{*}} \right),
\end{equation*}
where  $d_{\lambda}^*:=(-1)^{\frac{d_\lambda-1}{2}}\prod_{i=1}^k \lambda_i$.
Moreover, every $A_n$-character $\chi$ which is not algebraically conjugate to $\chi_{\gamma}$ has $\chi(\sigma) \in \Z.$
\end{enumerate}
\end{theorem}

\subsection{Some facts about partitions into distinct parts}

\begin{comment}
\begin{conjecture}[Erd\"os]
Let $p$ and $q$ be coprime integers.  Then there is a positive integer $N(p,q)$ such that every $n>N(p,q)$ may be expressed as a sum of distinct numbers of the form $p^a q^b$ over some set of distinct pairs of positive integers $a$ and $b$.
\end{conjecture}
\end{comment}
In 1959, Birch \cite{B} proved a conjecture of Erd\"os on representations of sufficiently large integers by sums of distinct numbers of the form $p^{a} q^{b}$, where $p,q$ are coprime and $a,b$ are positive integers.  In 2017, this result was quantified by Fang and Chen (see Theorem 1.1 of \cite{FC}).
\begin{theorem}[Fang-Chen]\label{Thm2Point2}
For any coprime integers $p,q>1$, there exist positive integers $K$ and $B$ with $$K<2^{2^{q^{2p}}},\hspace{.5cm}B<2^{2^{2^{q^{2p}}}}$$ such that every integer $n\geq B$ can be expressed as the sum of distinct terms taken from $$\left\{p^aq^b\mid a\geq0,\,0\leq b\leq K,\,a+b>0,\,a,b\in\mathbb{Z}\right\}.$$
\end{theorem}

Motivated by Birch's earlier work, Cassels (see Theorem 1 of  \cite{C}) studied the problem of representing sufficiently large integers as sums of distinct elements from
suitable integer sequences.
A careful inspection of his paper gives the following theorem, which was proved using the ``circle method".  
\begin{theorem}[Cassels]\label{Thm2Point3}
Suppose that $\mathbb{N}_T=\left\{a_1,a_2,\dots\right\}$ is a subset of positive integers in increasing order.  Let $S_T(n)$ denote the number of integers in $\mathbb{N}_T$ not exceeding $n$.  Choose $B_T\geq2^{2^{40}}$ to be a positive real number for which $$S_T(2n)-S_T(n)>2^{20}\log_2\log_2n$$ for all $n\geq B_T$, and choose $C_T\geq B_T$ to be a positive real number for which $$\sum_{s=1}^{S_T\left(C_T\right)}\left|\left|a_s\theta\right|\right|^2>2\log_2B_T+50$$ for all $\frac{1}{4B_T}\leq\theta\leq1-\frac{1}{4B_T}$.  If such $B_{T}$ and $C_{T}$ exist, then every positive integer $n>\max\left\{C_T^{10},2^{2^{2^{400}}}\right\}$ is expressible as a sum of distinct elements of $\mathbb{N}_T$.
\end{theorem}

\subsection{Proof of Theorem~\ref{Thm1}}

Let $\pi:=\left\{p_1,\dots,p_t\right\}$ be a set of distinct odd primes in increasing order.  For $n> N_{\pi}$, we will establish, for each $1\leq i\leq t$, that
$n-p_i$ is a sum of distinct squares of $\pi$-numbers. These representations imply that $n=p_i + M_i$, where $M_i$ is a sum of distinct
squares of $\pi$-numbers. Theorem~\ref{Thm2Point1}  then implies the theorem.

%Let $\lambda\left(p_i\right)\vdash n$ be the partition of $n$ represented by the above sum, $$\lambda\left(p_i\right):=\left(p_i,p_1^{2a_1},\dots,p_t^{2a_t}\right)$$ (with parts in non-increasing order), and let $\sigma_i\in A_n$ be a permutation with cycle type given by $\lambda\left(p_i\right)$.  Then by Theorem 2.1, we have that $$\mathbb{Q}\left(\chi\left(\sigma_i\right):1\leq i\leq t,\,\chi\in\mathrm{Irr}\left(A_n\right)\right)=\mathbb{Q}\left(\sqrt{d_{\lambda\left(p_i\right)}^*}:1\leq i\leq t\right).$$  Since all of the parts of each $\lambda\left(p_i\right)$ are perfect squares except for the part $p_i$ itself, it is clear that $$\mathbb{Q}\left(\sqrt{d_{\lambda\left(p_i\right)}^*}:1\leq i\leq t\right)=\mathbb{Q}\left(\sqrt{p_1^*},\sqrt{p_2^*},\dots,\sqrt{p_t^*}\right).$$

By Theorem~\ref{Thm2Point2}, for each pair $1\leq i < j\leq t$, we have that every integer $n\geq N_{i,j}:=2^{2^{2^{\Omega_\pi(i,j)}}},$ where $\Omega_\pi(i,j):=p_j^{4p_i^2}$, is a sum of distinct terms taken from the set $$\left\{p_i^{2a}p_j^{2b}:a,b\geq0,\,a+b>0,\,a,b\in\mathbb{Z}\right\}.$$
Obviously, the maximum bound occurs with $p_{t-1}$ and $p_{t}$, and so the theorem follows for $n> \mathcal{N}_{\pi}^{+}$. 
% In particular, for each $1\leq i\leq t$ we have that $$n=p_i+\sum_{\text{some }a_1,a_2,\dots,a_t\in\mathbb{Z}}p_1^{2a_1}p_2^{2a_2}\dots p_t^{2a_t}$$ for all $n\geq\mathcal{N}_{p_i}^+,$ where
%\begin{align*}
%\mathcal{N}_{p_i}^+&:=p_i+\max\left\{N_{i,j}:1 \leq i\neq j \leq t \right\}\\
%&=p_i+N_{t-1,t}.
%\end{align*}
%Choosing $N_{t-1,t}:=2^{2^{2^{\Omega_\pi}}},$ where $\Omega_\pi:=p_t^{4\left(p_{t-1}\right)^2},$ and letting $$\mathcal{N}_\pi^+:=p_t+N_{t-1,t},$$ we have that Theorem~\ref{Thm1} holds for all $n\geq\mathcal{N}_\pi^+$.
Finally, Cassels' result, where $\mathbb{N}_T$ denotes the set of squares of $\pi$-numbers,
guarantees that Theorem~\ref{Thm1} holds for all $n>\mathcal{N}_\pi^-$.  
%Finally, choosing $N_\pi:=\min\left\{\mathcal{N}_\pi^+,\mathcal{N}_\pi^-\right\}$ completes the proof.

\end{document}